\documentclass[12pt]{article}
\usepackage{amsfonts}
\usepackage{a4wide}
\usepackage{amsmath,amscd,amsthm,a4,amssymb}

\begin{document}

\begin{center}
{\large \textbf{A note on the representation of continuous \\
functions by linear superpositions}}

\

\textbf{Vugar E. Ismailov} \

{Mathematics and Mechanics Institute}

Azerbaijan National Academy of Sciences

Az-1141, Baku, Azerbaijan

{e-mail:} {vugaris@mail.ru}
\end{center}

\bigskip

\textbf{Abstract.} We consider the problem of the representation of real
continuous functions by linear superpositions $\sum_{i=1}^{k}g_{i}\circ
p_{i} $ with continuous $g_{i}$ and $p_{i}$. This problem was considered by
many authors. But complete, and at the same time explicit and practical
solutions to the problem was given only for the case $k=2$. For $k>2$, a
rather practical sufficient condition for the representation can be found in
Sternfeld [17] and Sproston, Strauss [16]. In this short note, we give a
necessary condition of such kind for the representability of continuous
functions.

\bigskip



\bigskip

\begin{center}
{\large \textbf{1. Introduction}}
\end{center}

Let $X$ be a set and $p_{i}$, $i=1,\ldots ,k$, be arbitrarily fixed
functions over $X$. For a given set $Y$, let $T(Y)$, $B(Y)$ and $C(Y)$ stand
for the space of all, bounded, and continuous real functions on $Y$,
respectively. Consider the following three sets

\begin{equation*}
S(X)=S(p_{1},\ldots,p_{k};X)=\left\{\sum_{i=1}^{k}g_{i}(p_{i}(x)):~x\in X,
~g_{i}\in T(\mathbb{R}),~i=1,\ldots,k\right\},
\end{equation*}

\begin{equation*}
S_{b}(X)=S_{b}(p_{1},\ldots,p_{k};X)=\left\{\sum_{i=1}^{k}g_{i}(p_{i}(x)):~x%
\in X, ~g_{i}\in B(\mathbb{R}),~i=1,\ldots,k\right\},
\end{equation*}

\begin{equation*}
S_{c}(X)=S_{c}(p_{1},\ldots,p_{k};X)=\left\{\sum_{i=1}^{k}g_{i}(p_{i}(x)):~x%
\in X, ~g_{i}\in C(\mathbb{R}),~i=1,\ldots,k\right\}.
\end{equation*}

For the second set, the functions $p_{i},i=1,\ldots ,k$, are considered to
be bounded on $X$. For the third set, we assume that $X$ is a compact
Hausdorff space and the functions $p_{i}$ are continuous on $X$. Members of
these sets will be called linear superpositions (see [20]).

At present, there are many works investigating possibilities of the
equalities $S(X)=T(X)$, $S_b(X)=B(X)$, $S_c(X)=C(X)$ (see [9] and references
therein). Here, we are interested in the last equality $S_c(X)=C(X)$.

The famous Kolmogorov superposition theorem states that for $X$ being the
unit cube in $\mathbb{R}^{d}$ there exist functions $p_{i}\in
C(X),~i=1,\ldots ,2d+1$, such that $S_{c}(p_{1},\ldots ,p_{2d+1};X)=C(X)$
(see [11]). Further the functions $p_{i}$ can be chosen as sums of
univariate functions. This deep result, which solved Hilbert's 13-th
problem, was generalized in many directions. One such direction was in
choosing various sets $X$ of $\mathbb{R}^{d}$, or even general metric spaces
(see, e.g., [3,4,9,14]). In all of these works, the functions $p_{1},\ldots
,p_{k}$ guaranteeing the equality $S_{c}(p_{1},\ldots ,p_{k};X)=C(X)$ were
incalculable. Regarding the nature of these functions, for some sets $X$,
they can be chosen to be at most from the class $Lip~1$ (see [5]).

Appearing in the late 60's, the work of Vitushkin and Henkin [20] showed
that for $p_{1},\ldots,p_{k}$ ($k$ may be very large) having except
continuity also smoothness properties, even the density of $%
S_{c}(p_{1},\ldots,p_{k};X)$ in $C(X)$ does not generally hold. Thus, the
question about when $S_{c}(X)=C(X)$ and $\overline{S_{c}(X)}=C(X)$ was
raised. Clearly, any answer depends on both the behavior of $%
p_{1},\ldots,p_{k}$ and the structure of $X$.

For the above problem of representation, the first crucial step was made by
Sternfeld [17]. He showed that the problem with its nature is dual to the
problem of uniform separation of measures of some certain class (see
[17,19]). The duality approach enabled him to prove that the number of terms
in the Kolmogorov superposition formula cannot be reduced (see [18]). Let $S$
be a class of measures defined on some field of subsets of $X$ and $F=\{p\}$
be a family of functions defined on $X$. $F$ uniformly separates measures of
the class $S$ if there exists a number $0<\lambda \leq 1$ such that for each
$\mu$ in $S$ the equality $\|\mu\circ p^{-1}\| \geq\lambda \|\mu\|$ holds
for some $p\in F$. In this terminology, $S_{b}(p_{1},\ldots,p_{k};X)=B(X)$
and $S_{c}(p_{1},\ldots,p_{k};X)=C(X)$ if and only if the family $%
\{p_1,\ldots,p_k\}$ uniformly separates measures of the classes $l_{1}(X)$
and $C(X)^{\ast}$ correspondingly (see [19]). Since $l_{1}(X)\subset
C(X)^{\ast}$ (the set of all regular Borel measures includes, in particular,
discrete measures), Sternfeld concluded that the equality $S_{c}(X)=C(X)$
implies $S_{b}(X)=B(X)$. In [7], we showed that any of these two equalities
implies $S(X)=T(X)$. That is, if some representation by linear
superpositions holds for continuous (or bounded) functions, then it holds
for all functions.

Sproston and Straus [16] gave a practically convenient sufficient condition
for the space $S_{c}(X)$ to be the whole of $C(X)$ (in fact, their result
was equivalently formulated in terms of sums of closed algebras). To
describe the condition, define the set functions
\begin{equation*}
\tau _{i}(Z)=\{x\in Z:~|p_{i}^{-1}(p_{i}(x))\bigcap Z|\geq 2\},\quad
Z\subset X,~i=1,\ldots ,k,
\end{equation*}%
where $|Y|$ denotes the cardinality of a considered set $Y$. Define $\tau
(Z) $ to be $\bigcap_{i=1}^{k}\tau _{i}(Z)$ and define $\tau ^{2}(Z)=\tau
(\tau (Z))$, $\tau ^{3}(Z)=\tau (\tau ^{2}(Z))$ and so on inductively. The
result of [16] says that $S_{c}(p_{1},\ldots ,p_{k};X)=C(X)$ provided that $%
\tau ^{n}(X)=\emptyset $ for some positive integer $n$. In fact, this
condition first appeared in the work of Sternfeld [17], where the author
proved that $\tau ^{n}(X)=\emptyset $ (for some $n$) guarantees that the
family $\{p_{1},\ldots ,p_{k}\}$ uniformly separates measures of the class $%
l_{1}(X)$ and also regular Borel measures if $X$ is a compact metric space.
Sproston and Straus proved the last statement for $X$ being a compact
Hausdorff space. For $k=2$, the condition is also necessary for the
representation, but not in general if $k>2$ (see the counterexample in [16]).

For $k=2$, the above condition $\tau ^{n}(X)=\emptyset $ can be
expressed in terms of sets of points in $X$ which were introduced
in the literature under different names such as ``bolts of
lightning'' [8,9,13], ``trips'' [12], ``paths'' [6], ``loops''
[2], etc. These objects are geometrically explicit. A path with
respect to two functions $p_{1}$ and $p_{2}$ can be described as a
trace of some point in $X$ traveling (more precisely, jumping) in
alternating level sets of the functions $p_{1}$ and $p_{2}$. If
the point returns to its primary position after such a travel, the
obtained set is called a closed path. It is not difficult to prove
that $\tau ^{n}(X)=\emptyset $ if and only if there are no closed
paths in $X$ and the lengths (number of points) of all paths are
uniformly bounded (see [9]).

Paths with respect to two coordinate functions (and two algebras) have been
extensively implemented by Marshall and O'Farrell [12,13] to solve the
problem concerning density of $S_{c}(p_{1},p_{2};X)$ in $C(X)$. Their work
[13] showed the essence of such paths by explaining that every regular Borel
measure orthogonal to $S_{c}(p_{1},p_{2};X)$ in $C(X)$ is in the closure of
the set of measures generated by paths.

It should be remarked that many authors considering the problems of
representation and approximation by linear superpositions indicated the
difficulties and at the same time usefulness of going from measure-theoretic
to path-descriptive results (see, e.g., [9,13,19]).

The purpose of this note is to obtain a path-descriptive necessary condition
for representability of each continuous function by linear superpositions.
We hope that our condition will complement the above sufficient condition $%
\tau^{n}(X)=\emptyset$ in some sense.

\bigskip

\begin{center}
{\large \textbf{2. The main result}}
\end{center}

We begin this section with a definition. Let us first assume that we are
given a compact Hausdorff space $X$ and continuous functions $%
p_i:X\rightarrow \mathbb{R}$, $i=1,\ldots,k$.

\bigskip

\textbf{Definition 2.1} (see [1,7,10]). \textit{A set of points $%
l=(x_{1},\ldots ,x_{n})\subset X$ is called a closed path with respect to
the functions $p_{1},\ldots ,p_{k}$ if there exists a vector $\lambda
=(\lambda _{1},\ldots ,\lambda _{n})\in \mathbb{Z}^{n}\setminus \{\mathbf{0}%
\}$ such that
\begin{equation*}
\sum_{j=1}^{n}\lambda _{j}\delta _{p_{i}(x_{j})}=0,\quad for~all~i=1,\ldots
,k.
\end{equation*}%
}

Here $\delta _{a}$ is the characteristic function of the set $\{a\}$.

For example, the set $l=\{(0,0,0),~(0,0,1),~(0,1,0),~(1,0,0),~(1,1,1)\}$ is
a closed path in $\mathbb{R}^{3}$ with respect to the functions $%
p_{i}(z_{1},z_{2},z_{3})=z_{i},~i=1,2,3.$ The vector $\lambda $ in
Definition 2.1 can be taken as $(-2,1,1,1,-1).$

The idea of closed paths with respect to $k$ directions in $\mathbb{R}^{d}$
was first considered in the paper by Braess and Pinkus [1]. Klopotowski,
Nadkarni, Rao [10] defined these objects with respect to canonical
projections. In our recent paper [7], which deals with linear
superpositions, closed paths have been generalized to those having
association with $k$ arbitrary functions. In these three works, it was shown
that nonexistence of closed paths of the respective form is both necessary
and sufficient for

1) interpolation by ridge functions [1];

2) representation of multivariate functions by sums of univariate functions
[10];

3) representation by linear superpositions [7].

It should be remarked that consideration of only closed paths is not enough
for investigating the problems of representation by linear superpositions in
cases when some topology (that of boundedness, or continuity) is involved.
As in the case $k=2$, more general objects must be implemented.

\bigskip

\textbf{Definition 2.2.} \textit{A set of points $l=(x_{1},\ldots
,x_{n})\subset X$ is called a path with respect to the functions $%
p_{1},\ldots ,p_{k}$ if there exists a vector $\lambda =(\lambda _{1},\ldots
,\lambda _{n})\in \mathbb{Z}^{n}\setminus \{\mathbf{0}\}$ such that for any $%
i=1,\ldots ,k$
\begin{equation*}
\sum_{j=1}^{n}\lambda _{j}\delta _{p_{i}(x_{j})}=\sum_{s=1}^{r_{i}}\lambda
_{i_{s}}\delta _{p_{i}(x_{i_{s}})},\quad where~r_{i}\leq k.
\end{equation*}%
}

Note that for $i=1,\ldots ,k$, the set $\{\lambda _{i_{s}},~s=1,...,r_{i}\}$
is a subset of the set $\{\lambda _{j},~j=1,...,n\}$. Thus, for each $i$, we
actually have at most $k$ terms in the sum $\sum_{j=1}^{n}\lambda _{j}\delta
_{p_{i}(x_{j})}$.

Let, for example, $k=2,$ $p_{1}(x_{1})=p_{1}(x_{2})$, $
p_{2}(x_{2})=p_{2}(x_{3})$, $ p_{1}(x_{3})=p_{1}(x_{4})$,...,
$p_{2}(x_{n-1})=p_{2}(x_{n})$. Then it is not difficult to see
that for a vector $\lambda =(\lambda _{1},\ldots ,\lambda _{n})$
with the components $\lambda _{i}=(-1)^{i},$
\begin{eqnarray*}
\sum_{j=1}^{n}\lambda _{j}\delta _{p_{1}(x_{j})}=\lambda
_{n}\delta _{p_{1}(x_{n})}, \\
\sum_{j=1}^{n}\lambda _{j}\delta
_{p_{2}(x_{j})}=\lambda _{1}\delta _{p_{2}(x_{1})}.
\end{eqnarray*}
Thus, by Definition 2.2, the set $l=\{x_{1},\ldots ,x_{n}\}$ forms a path
with respect to the functions $p_{1}$ and $p_{2}$.

One can construct many paths by adding not more than $k$ arbitrary points to
a closed path with respect to some functions $p_{1},\ldots ,p_{k}.$

\bigskip

\textit{Remark 1}. Closed paths and paths with respect to two functions
mentioned above in the Introduction satisfy Definitions 2.1 and 2.2,
respectively. But for $k=2$, Definitions 2.1 and 2.2 may allow also some
unions of the previously known objects.

\bigskip

Each path $l=(x_1,\ldots,x_n)$ and the associated vector $%
\lambda=(\lambda_1,\ldots,\lambda_n)$ generate the functional

\begin{equation}
G_{l,\lambda}(f)=\sum_{j=1}^{n}\lambda_{j}f(x_j), \quad f\in C(X).
\end{equation}

Clearly, $G_{l,\lambda }$ is linear and continuous with norm $%
\sum_{j=1}^{n}|\lambda _{j}|.$

Note that the space $S_{c}(p_{1},\ldots ,p_{k};X)$ is the sum of algebras
\begin{equation*}
S_{i}=\{u_{i}\in C(X):~u_{i}=g(p_{i}(x)),~g\in C(\mathbb{R})\},\quad
i=1,\ldots ,k.
\end{equation*}%
From Definition 2.2 it follows that for each function $u_{i}\in S_{i}$, $%
i=1,\ldots ,k$,
\begin{equation}
G_{l,\lambda }(u_{i})=\sum_{j=1}^{n}\lambda
_{j}u_{i}(x_{j})=\sum_{s=1}^{r_{i}}\lambda _{i_{s}}u_{i}(x_{i_{s}}),
\end{equation}%
where $r_{i}\leq k$. That is, for each algebra $S_{i}$, $G_{l,\lambda }$ can
be reduced to a functional defined with the help of not more than $k$ points
of the path $l$. Note that if $l$ is closed, then $G_{l,\lambda }(u_{i})=0$
for all $u_{i}\in S_{i}$, $i=1,\ldots ,k$, whence $G_{l,\lambda }(u)=0$, for
any $u\in S_{c}(p_{1},\ldots ,p_{k};X)$.

\bigskip

\textit{Remark 2.} Let $f\in C(X)$. If $G_{l,\lambda}(f)=0$, for any closed
path $l\subset X$, then $f\in S(p_{1},\ldots,p_{k};X)$. That is, $f$ can be
represented by linear superpositions. But generally, $f$ may not be in $%
S_{c}(p_{1},\ldots,p_{k};X)$ (see [7]).

\bigskip

\textbf{Theorem 2.3.} \textit{Let $S_{c}(p_{1},\ldots,p_{k};X)=C(X)$. Then}

\textit{(a) there are no closed paths in $X$.}

\textit{(b) lengths (number of points) of all paths in $X$ are uniformly
bounded.}

\begin{proof}

The part (a) is obvious. Indeed, let $l=(x_1,\ldots,x_n)$ be a closed path
in $X$ and $\lambda=(\lambda_1,\ldots,\lambda_n)$ be a vector associated
with it. As it is indicated above, $G_{l,\lambda}(u)=0$, for any function $%
u\in S_{c}(X)$. Let $f_0$ be a continuous function such that $f_0(x_j)=1$ if
$\lambda_j>0$ and $f_0(x_j)=-1$ if $\lambda_j<0$, $j=1,\ldots,n$. Since $%
G_{l,\lambda}(f_0)\neq 0$, $f_0$ cannot be in $S_{c}(X)$. Therefore, $%
S_{c}(X)\neq C(X)$. But this contradicts the hypothesis of the theorem.

To prove the (b)-part of the assertion, consider the linear space
\begin{equation*}
U=\prod_{i=1}^{k}S_i=\{(u_1,\ldots,u_k): ~u_i\in S_i, ~i=1,\ldots,k\}
\end{equation*}
endowed with the norm
\begin{equation*}
\|(u_1,\ldots,u_k)\|=\|u_1\|+\cdots+\|u_k\|.
\end{equation*}

We will also deal with the dual of the space $U$. Each functional $F\in
U^{\ast}$ can be written as the sum
\begin{equation*}
F=F_1+\cdots+F_k,
\end{equation*}
where the functionals $F_i\in S_i^{\ast}$ and
\begin{equation*}
F_i(u_i)=F[(0,\ldots,u_i,\ldots,0)], \quad i=1,\ldots,k.
\end{equation*}
Thus, we see that the functional $F$ determines the collection $%
(F_1,\ldots,F_k)$. Conversely, every collection $(F_1,\ldots,F_k)$ of
continuous linear functionals $F_i\in S_i^{\ast}$, $i=1,\ldots,k$,
determines the functional $F_1+\cdots+F_k,$ on $U$. Considering this, in
what follows, elements of $U^{\ast}$ will be denoted by $(F_1,\ldots,F_k)$.

It is not difficult to verify that
\begin{equation}
\|(F_1,\ldots,F_k)\|=max\{\|F_1\|,\ldots,\|F_k\|\}.
\end{equation}

Consider the operator
\begin{equation*}
A:U\rightarrow C(X),\quad A[(u_{1},\ldots ,u_{k})]=u_{1}+\cdots +u_{k}.
\end{equation*}%
Clearly, $A$ is a linear continuous operator with the norm $\Vert A\Vert =1$%
. Besides, since $S_{c}(X)=C(X)$, $A$ is a surjection. Consider also the
conjugate operator
\begin{equation*}
A^{\ast }:C(X)^{\ast }\rightarrow U^{\ast },~A^{\ast }[H]=(F_{1},\ldots
,F_{k}),
\end{equation*}%
where $F_{i}(u_{i})=H(u_{i})$, for any $u_{i}\in S_{i}$, $i=1,\ldots ,k$.
Let H be an arbitrary functional $G_{l,\lambda }$ of the form (1). Set $%
A^{\ast }[G_{l,\lambda }]=(G_{1},\ldots ,G_{k})$. From (2) it follows that
\begin{equation*}
|G_{i}(u_{i})|=|G_{l,\lambda }(u_{i})|\leq \Vert u_{i}\Vert \times
\sum_{s=1}^{r_{i}}|\lambda _{i_{s}}|\leq b_{\lambda }(k)\times \Vert
u_{i}\Vert ,\quad i=1,\ldots ,k,
\end{equation*}%
where $b_{\lambda }(k)$ is the maximum of all numbers $\sum_{s=1}^{k}|%
\lambda _{j_{s}}|$ formed by $k$ components of the vector $\lambda $.
Therefore,
\begin{equation*}
\Vert G_{i}\Vert \leq b_{\lambda }(k),\quad i=1,\ldots ,k.
\end{equation*}%
From (3) we obtain that
\begin{equation}
\Vert A^{\ast }[G_{l,\lambda }]\Vert =\Vert (G_{1},\ldots ,G_{k})\Vert \leq
b_{\lambda }(k)
\end{equation}%
Since $A$ is a surjection, there exists a positive real number $\delta $
such that
\begin{equation*}
\Vert A^{\ast }[H]\Vert >\delta \Vert H\Vert
\end{equation*}%
for any functional $H\in C(X)^{\ast }$ (see Rudin [15]). Taking into account
that $\Vert G_{l,\lambda }\Vert =\sum_{j=1}^{n}|\lambda _{j}|$, for the
functional $G_{l,\lambda }$ we have
\begin{equation}
\Vert A^{\ast }[G_{l,\lambda }]\Vert >\delta \sum_{j=1}^{n}|\lambda _{j}|.
\end{equation}%
It follows from (4) and (5) that
\begin{equation*}
\delta <\frac{b_{\lambda }(k)}{\sum_{j=1}^{n}|\lambda _{j}|}.
\end{equation*}%
The last inequality shows that $n$ (the length of the arbitrarily
chosen path $l$) cannot be as great as possible, otherwise $\delta
=0$. This simply means that there must be some positive integer
bounding the lengths of all paths in $X$. \end{proof}

\bigskip

\textit{Remark 3.} The condition (a) of Theorem 2.3 is also necessary for
the density of $S_{c}(p_{1},\ldots,p_{k};X)$ in $C(X)$, whereas the
condition (b) is not. For the case $k=2$, one can use the nontrivial example
of Khavinson [8].

\bigskip

\textit{Remark 4.} The conditions (a) and (b) of Theorem 2.3 are sufficient
for the equality $S_{c}(p_{1},p_{2};X)=C(X)$ (see [9]). Thus in the case $%
k=2 $, these conditions are equivalent to the above mentioned condition $%
\tau ^{n}(X)=\emptyset $ of Sternfeld. For $k>2$, they are no longer
equivalent, since the condition of Sternfeld is not necessary for the
representation (see [16]), The question if for $k>2$, our conditions (a) and
(b) are sufficient for the representation, unfortunately has a negative
answer. Our argument is as follows. It can be proven by the same way that
the conditions (a) and (b) are necessary for the equality $%
S_{b}(p_{1},\ldots ,p_{k};X)=B(X) $. If they had been sufficient for $%
S_{c}(p_{1},\ldots ,p_{k};X)=C(X)$, they would have been also sufficient for
$S_{b}(p_{1},\ldots ,p_{k};X)=B(X)$, since the representability of
continuous functions implies the representability of bounded functions (see
[17]). Hence, we would obtain that the conditions (a) and (b) are necessary
and sufficient for both the equalities $S_{c}(p_{1},\ldots ,p_{k};X)=C(X)$
and $S_{b}(p_{1},\ldots ,p_{k};X)=B(X)$. But for $k>2$, these equalities are
not equivalent (see [19]).

\bigskip

\textbf{Acknowledgments.} I would like to thank the anonymous
referee for his/her comments on improving the original manuscript.

\end{document}